\input amstex
\documentstyle{amsppt} 
\TagsOnRight 
\NoBlackBoxes

\def\Prob{\text{Prob }} 
\def\RR{\text{{\rm I}\!{\rm R}}}
\def\NN{\text{{\rm I}\!{\rm N}}}
\def\Rd{\RR^d} 
\def\CC{C} 
\def\t{\text}

\def\pa{\partial}
\def\ra{\rangle}
\def\la{\langle}

\def\loc{\,\text{\rm loc}\,}
\def\sgn{\,\text{\rm sgn}\,}
\def\supp{\,\text{\rm supp}\,}

\def\div{\operatorname{div}}
\def\grad{\nabla}
\def\eps{\varepsilon}
\def\del{\delta}
\def\lap{\Delta}

\define\refCLF           {1}
\define\refDiPerna       {2}
\define\refHLone         {3}
\define\refHLtwo         {4}
\define\refKruzkov       {5}
\define\refLL            {6}
\define\refLN            {7}
\define\refVonNeumann    {8}
\define\refSchonbek      {9}
\define\refSzepessy      {10}
\define\refTartar        {11}
\define\refVolpert       {12} 

\topmatter
\title 
Nonlinear diffusive-dispersive limits                         
for multidimensional conservation laws   
\endtitle
\footnote""{
Published in: \bf  ``Advances in nonlinear partial differential equations and related areas''
 (Beijing, 1997), pp.~103--123, 
World Sci. Publ., River Edge, NJ, 1998. 
\hskip.2cm \hfill}
\author Joaquim M. Correia
\footnote "$^1$"
{ \, 
Former affiliation: 
Ecole Polytechnique, Palaiseau. 
\hskip8.cm \hfill}
and Philippe G. LeFloch${^{1,}}$
\footnote"$^{2}$"
{ \, Current affiliation (from 2004): {\it Laboratoire Jacques-Louis Lions
        \& Centre National de la Recherche Scientifique (CNRS), 
        Universit\'e Pierre et Marie Curie (Paris 6), 4 Place Jussieu, 
        75252 Paris, France.}  E-mail: { pgLeFloch\@gmail.com}.
\hskip1.cm \hfill
\hskip.1cm \hfill}
\endauthor
\rightheadtext{Nonlinear diffusive-dispersive limits}
\leftheadtext{}
\abstract 
We consider a class of multidimensional conservation
laws with vanishing nonlinear diffusion and dispersion terms. 
Under a condition on the relative size of the diffusion and dispersion 
coefficients, we establish that the diffusive-dispersive solutions 
are uniformly bounded in a space $L^p$ ($p$ arbitrary large, depending 
on the nonlinearity of the diffusion) and 
converge to the classical, entropy solution of the corresponding 
multidimensional, hyperbolic conservation law. Previous results were 
restricted to one-dimensional equations and specific spaces $L^p$. 
Our proof is based on DiPerna's 
uniqueness theorem in the class of entropy measure-valued solutions. 
\endabstract 
\endtopmatter 

\document 

\

\subheading{1. Introduction}

Nonlinear hyperbolic conservation laws arise in the modeling of many problems 
from continuum mechanics, physics, chemistry, etc. The equations become 
parabolic when additional dissipation mechanisms are taken into account: 
diffusion, heat conduction, 
capillarity in fluids, Hall effect in magnetohydrodynamics, etc. From 
a mathematical standpoint, hyperbolic equations admit discontinuous 
solutions while parabolic equations have smooth solutions. 
Discontinuous solutions, understood in the generalized sense of the 
distribution theory, are usually non-unique. It is therefore fundamental 
to understand which solutions are selected by a specific zero
diffusion-dispersion 
limit. In this paper we address this issue for multidimensional, 
scalar conservation laws, and review previous work on the subject restricted to 
one-dimensional equations. 

Consider the Cauchy problem
$$
\gather
\pa_t u + \div\!f(u) = 0, \qquad (x,t)\in\Rd\times\RR_+,
\tag 1.1a \\
u(x,0) = u_0(x), \qquad x\in\Rd,
\tag 1.1b
\endgather
$$
where the unknown function $u$ is scalar-valued. Smooth solutions to (1.1) 
also satisfy an infinite list of additional conservation laws~: 
$$
\pa_t \eta(u) + \div q(u) = 0, \qquad \grad_u q = \grad_u \eta \, \grad_u\!f,
\tag 1.2
$$
where $\eta$ is a convex function of $u$. For discontinuous solutions, 
Kru\v{z}kov \cite{\refKruzkov} shows that (1.2) should be replaced 
by the set of inequalities 
$$
\pa_t \eta(u) + \div q(u) \leq 0,
\tag 1.3
$$
which select physically meaningful, discontinuous solutions. 
The condition (1.3) is called an entropy inequality; 
it is motivated by the second law of thermodynamics, in the context of 
gas dynamics. By definition, 
an entropy solution of problem (1.1) satisfies (1.1) in the sense of 
distributions, and additionally  
(1.3) for any entropy pair $(\eta,q)$ with convex $\eta$.

Consider the following approximation of (1.1) obtained by 
adding a nonlinear diffusion, $b:~\Rd \to \Rd$, 
and a linear dispersion to the right hand side of (1.1a)~: 
$$
\gather
\pa_t u + \div\!f(u) = \div \left(
   \eps\,b_j(\grad u) + \del\,\pa^2_{x_j} \! u \right)_{1 \leq j \leq d}, 
	  \qquad (x,t)\in\Rd\times\RR_+,
\tag 1.4a \\
u(x,0) = u_0^{\eps,\del}(x),
   \qquad x\in\Rd. 
\tag 1.4b
\endgather
$$
Let $u^{\eps,\del}:\Rd\times[0,T]\to\RR$ be smooth solutions
defined on an interval $[0,T]$ with a uniform $T$ independent of $\eps, \del$. 
In (1.4b), 
$u_0^{\eps,\del}$ is an approximation of the initial condition $u_0$ in (1.1b).

Our main purpose is to derive conditions under which, as $\eps>0$ and $\del$
tend
to zero,
the solutions $u^{\eps,\del}$ converge in a strong topology to the entropy
solution
of (1.1). When $\eps = 0$, equation (1.4a) is a generalized version of 
the well-known Korteweg-de Vries (KdV) equation, and the solutions become more
and 
more oscillatory as $\del~~\to~~0$~: the approximate solutions do not converge 
to zero; see Lax and Levermore \cite{\refLL}. When 
$\del = 0$, (1.4a) reduces to a nonlinear parabolic equation, resembling the 
pseudo-viscosity approximation of von Neumann and Richtmyer
\cite{\refVonNeumann}; in that regime, the solution converges strongly
to the entropy solution. Therefore, 
to ensure the convergence of the zero diffusion-dispersion approximation (1.4), 
it is necessary that diffusion dominate dispersion. Indeed the main result 
of the present paper establishes, under rather broad assumptions 
(see Section 3, Theorem 3.1--3.3), that 
the solution of (1.4) tends to the classical entropy solution of (1.1) 
when $\eps, \delta \to 0$ with $|\delta| << \eps$. 

For clarity, the main assumptions made in this paper are collected here. 
First concerning the flux function we shall assume 
\roster
\item"{ ($H_1$) }" for some $C_1, C_1'>0$ and $m\ge0$, \quad
                       $\left| f'(u) \right| \leq
                       C_1 +C_1'\left| u \right|^{m-1}$
							  \quad for all $u\in\RR$.
\endroster
For the diffusion term, we fix $r\ge0$ and assume 
\roster
\item"{ ($H_2$) }" for some $C_2>0$, \quad
                       $C_2 \left| \lambda \right|^{r+1} \leq
                       \lambda \cdot\,b(\lambda) \leq 
							  C_3 \left| \lambda \right|^{r+1}$
							  \quad for all $\lambda\in\Rd$. 
\endroster
In the case $0\leq r<2$, we will need also 
\roster
\item"{ ($H_3$) }" $Db(\lambda)$ is a positive definite 
                   matrix uniformly in $\lambda\in\Rd$.
\endroster
We remark that the diffusion $b_j(\grad u) = \partial_{x_j} u$ 
satisfies ($H_3$).

The case $d=1$ of one-dimensional equations was treated in the important paper 
by Schonbek \cite{\refSchonbek}, where, in particular, the concept of 
$L^p$ Young measures is introduced together with an extension of the 
compensated compactness method for conservation laws. We follow here
LeFloch and Natalini \cite{\refLN} who, also for
one-dimensional 
equations, developed another approach based on DiPerna's uniqueness theorem 
for entropy measure-valued solutions \cite{\refDiPerna} 
(see Section 2 for a review). 
Specifically one uses a generalization of DiPerna's result to $L^p$ functions, 
due to Szepessy \cite{\refSzepessy}. The present paper therefore relies on a 
method of proof that was successful first in proving convergence 
of finite difference schemes~: Szepessy (\cite{\refSzepessy} 
and the references therein by Szepessy and co-authors) 
and Coquel and LeFloch \cite{\refCLF}. 

Recent work by Hayes and LeFloch (see \cite{\refHLone,\refHLtwo}) 
treats the transitional case where both terms, the diffusion and the
dispersion, are in balance. Convergence results in this regime cannot be
obtained by 
the measure-valued solutions approach. 


\

\subheading{2. Entropy Measure-Valued Solutions} 

We include here, as background for further reference,
basic material on Young measures and entropy measure-valued
(e.m.-v.) solutions. First of all we will need 
Schonbek's representation theorem
for the Young measures associated with a sequence of
uniformly bounded in $L^q$.
The corresponding setting in $L^\infty$ was first established
by Tartar \cite{\refTartar}. 
In the whole of this subsection, $q \in (1,\infty)$ and
$T\le\infty$ are fixed.

\proclaim{Lemma 2.1} {\rm (See \cite{\refSchonbek}.)}
Let $\{u_k\}$ be a bounded sequence in
$L^\infty((0,T); L^q (\Rd))$. Then there exists a subsequence
still denoted by $\{u_k\}$ and a weakly $\star$ measurable 
mapping $\nu:\Rd\times(0,T)\to \Prob(\RR)$ taking 
its values in the space of non-negative measures
with unit total mass (probability measures) such that,
for all functions $g\in \CC(\RR)$ satisfying 
$$
g(u) = \Cal O(|u|^s)  \quad \t{\rm as }|u|\to \infty,
       \qquad \t{ for some $s\in [0,q)$,} 
\tag 2.1
$$
the following limit representation holds
$$
\lim_{k\to\infty}\iint_{\Rd\times (0,T)} g(u_k(x,t))\,\phi(x,t)\,\,dxdt
=
\iint_{\Rd\times(0,T)} \la \nu_{(x,t)},g \ra \,\,\phi(x,t)\,\,dxdt 
\tag 2.2
$$
for all $\phi\in L^1(\Rd\times(0,T))\cap L^\infty(\Rd\times(0,T))$.

Conversely, given $\nu$, there exists a sequence $\{u_k\}$ 
satisfying the same conditions as above such that {\rm (2.2)} 
holds for any $g$ satisfying {\rm (2.1)}.
\endproclaim

We use the notation
$\la \nu_{(x,t)},g \ra := \int_{\RR} g(u)\,d\nu_{(x,t)}$,
which therefore describes $weak\star-\lim g(u_k)$. ``Weak $\star$
measurable'' means that the real-valued function
$\la \nu_{(x,t)},g \ra$ is measurable with respect to $(x,t)$,
for each continuous $g$ satisfying (2.1).
The measure-valued function $\nu_{(x,t)}$ is called a Young measure
associated with the sequence $\{u_k\}$. The following result
reveals the connection between the structure of $\nu$ and the
strong convergence of the subsequence.

\proclaim{Lemma 2.2} Suppose that $\nu$ is a Young measure
associated with a sequence $\{u_k\}$, bounded
in $L^\infty((0,T); L^q (\Rd))$. 
For $u \in L^\infty((0,T); L^q (\Rd))$, the following
statements are equivalent:
\roster
\item"{(i)}" $\lim_{k\to \infty}u_k = u$
      \quad \t{in $L^s((0,T); L^p_{\loc}(\Rd))$,}
		\qquad \t{ for all $s<\infty$ and $p\in [1,q)$;}
\item"{(ii)}" $\nu_{(x,t)} = \delta_{u (x,t)}$
      \qquad \t{ for a.e. $(x,t)\in\Rd\times(0,T)$.}
\endroster
\endproclaim

\noindent In (ii) above, the notation 
$\delta_{u(x,t)}$ 
is used for the Dirac mass defined by
$$
\la \delta_{u(x,t)}, g \ra = g(u(x,t)) \qquad
\text{ for all $g\in \CC(\RR)$ satisfying (2.1).}
$$
Following DiPerna \cite{\refDiPerna} and Szepessy \cite{\refSzepessy}, 
we define the e.m.-v. solutions to the
first order Cauchy problem (1.1).

\proclaim{Definition 2.1} 
Assume that $f\in\CC(\RR)^d$ satisfies the growth condition
{\rm (2.1)} and $u_0\in L^1(\Rd)\cap L^q (\Rd)$. 
A Young measure $\nu$ associated with a sequence $\{u_k\}$,
which is assumed to be bounded in $L^\infty( (0,T);L^q (\Rd))$,
is called an entropy measure-valued (e.m.-v.) solution
to the problem {\rm (1.1)} if
$$
\pa_t \la \nu_{(\cdot)},|u-k| \ra +
\div \la \nu_{(\cdot)},\sgn(u-k)(f(u)-f(k)) \ra \leq 0 
\qquad \t{ for all $k\in\RR$},
\tag 2.5a
$$
in the sense of distributions on $\Rd\times (0,T)$ and
$$
\lim_{t\to 0^+}{1\over t} \! \int^t_0 \!\! \int_K \la \nu_{(x,s)},
|u-u_0(x)| \,\ra \,dxds = 0,
\qquad \t{ for all compact set $K\subseteq \Rd$.}
\tag 2.5b
$$
\endproclaim

A function $u\in L^\infty( (0,T);L^1(\Rd)\cap L^q (\Rd))$ is 
an entropy weak solution to (1.1) in the sense of Kru\v{z}kov 
\cite{\refKruzkov} and Volpert \cite{\refVolpert} if and only if
the Dirac measure $\delta_{u(\cdot)}$ is an e.m.-v. solution.
In the case $q=+\infty$, existence and uniqueness of such solutions
was shown in \cite{\refKruzkov}. The following results on
e.m.-v. solutions were established in \cite{\refSzepessy}:
Proposition 2.3 states that e.m.-v. solutions are actually
Kru\v{z}kov's solutions. 
Proposition 2.4 states that the problem has a unique solution in $L^q$. 

\proclaim{Proposition 2.3} 
Assume that $f$ satisfies {\rm (2.1)} and
$u_0\in L^1(\Rd)\cap L^q (\Rd)$. 
Suppose that $\nu$ is an e.m.-v. solution to {\rm (1.1)}.
Then there exists a function
$u\in L^\infty( (0,T);L^1(\Rd)\cap L^q (\Rd))$ such that
$$
\nu_{(x,t)}=\delta_{u(x,t)} 
\qquad \t{ for a.e. $(x,t)\in\Rd\times (0,T)$.}
\tag 2.6
$$
\endproclaim

\proclaim{Proposition 2.4} 
Assume that $f$ satisfies {\rm (2.1)} 
and $u_0\in L^1(\Rd)\cap L^q (\Rd)$. 
Then there exists a unique entropy solution
$$
u \in L^\infty( (0,T);L^1(\Rd)\cap L^q (\Rd))
$$
to {\rm (1.1)} which, moreover, satisfies 
$$
\|u(t)\|_{L^p(\Rd)}\leq \|u_0\|_{L^p(\Rd)}
\qquad \t{ for a.e. $t\in\, (0,T)$ and all $p\in[1,q]$.}
\tag 2.7
$$
The measure-valued mapping $\nu_{(x,t)}=\delta_{u(x,t)}$ is 
the unique e.m.-v. solution of the same problem.
\endproclaim

Combining Propositions 2.3 and 2.4 and Lemma 2.2, we obtain the main 
convergence tool~:

\proclaim{Corollary 2.5} 
Assume that $f$ satisfies {\rm (2.1)} and
$u_0\in L^1(\Rd)\cap L^q (\Rd)$ for $q>1$. 
Let $\{u_k\}$ be a sequence bounded in
$L^\infty( (0,T);L^q (\Rd))$ and let $\nu$
be a Young measure associated with this sequence. 
If $\nu$ is an e.m.-v. solution to {\rm (1.1)}, then
$$
\lim_{k\to\infty}u_k = u
\quad \t{ in } L^s( (0,T);L^p_{loc}(\Rd)) 
\qquad \t{ for all $s<\infty$ and all $p\in [1,q)$,}
$$
where $u \in L^\infty( (0,T);L^1 (\Rd)\cap L^q (\Rd))$
is the unique entropy solution to {\rm (1.1)}.
\endproclaim


\

\subheading{3. Convergence Results}

Throughout it is assumed $u_0 \in L^1(\Rd) \cap L^q(\Rd)$ 
and the initial data in (1.4b) are smooth functions
with compact support and are uniformly bounded in $L^1(\Rd) \cap L^q(\Rd)$
for some $q>2$. While in previous works 
\cite{\refSchonbek, \refLN}, a single value of $q$ was treated, 
we can here handle arbitrary large values of $q$. For simplicity in the presentation, 
we will always consider exponents $q$ of the form 
$$
q=2+n(r-1),
$$
where $n \geq 0$ is any integer. 
Therefore, when the diffusion is superlinear, in the sense that ($H_2$) holds with 
$r>1$, then arbitrary large values of $q$ are obtained.  
Restricting attention to the diffusion-dominant 
regime $\del = O(\eps)$, we suppose that 
$u_0^{\eps,\del}$ approaches the initial condition $u_0$ of (1.1b) in the sense 
that~:
$$
\aligned
&\lim_{\eps \to 0+} u_0^{\eps,\del} = u_0
\quad \t{ in } L^1(\Rd) \cap L^q(\Rd), \\
&\| u_0^{\eps,\del} \|_{L^2(\Rd)} \leq \| u_0 \|_{L^2(\Rd)}.
\endaligned
\tag 3.1
$$
The following convergence theorems concern a sequence $u^{\eps,\del}$ 
of smooth solutions to problem {\rm (1.4)}, 
defined on $\Rd\times[0,T]$ and decaying rapidly at
infinity. 

First consider the hypothesis ($H_2$) with $r\geq 2$, that is the case of 
diffusions with (at least) quadratic growth. 

\proclaim{Theorem 3.1} Suppose that the flux $f$ satisfies {\rm ($H_1$)} with $m<q$
(which is always possible when $r>1$ by taking $q$ large enough). 
Suppose that 
the diffusion $b$ satisfies {\rm ($H_2$)} with $r\ge2$. 
If $\del = o(\eps^\frac{3}{r+1})$, then the sequence $u^{\eps,\del}$ 
converges in $L^s\left((0,T);L^p(\Rd) \right)$, for all $s<\infty$ and 
$p<q$, to a function
$$
u \in L^\infty\left(  (0,T);L^1(\Rd) \cap L^q(\Rd) \right), 
$$ 
which is the unique entropy solution to {\rm (1.1)}.
\endproclaim

Observe that $m$ and $q$ can be arbitary large in Theorem 3.1. 
To treat the case $r < 2$, we need the additional condition ($H_3$) 
on the diffusion. First for diffusion with linear growth ($r=1$), 
we obtain a result in the space $L^2$~:

\proclaim{Theorem 3.2} Suppose that $f$ satisfies {\rm ($H_1$)} with $m \leq 1$,
and $b$ satisfies {\rm ($H_2$)}-{\rm ($H_3$)} with $r = 1$.
If $\del = o(\eps^2)$, then the sequence $u^{\eps,\del}$ 
converges in $L^s\left((0,T);L^p(\Rd) \right)$, 
for all $s<\infty$ and $p<2$, to a function
$$
u\in L^\infty\left((0,T);L^1(\Rd) \cap L^2(\Rd) \right),
$$ 
which is the unique entropy solution to {\rm (1.1)}.
\endproclaim

In particular Theorem 3.2 covers the interesting case of a linear diffusion and 
a linear dispersion with an (at most) linear flux at infinity. 
More generally, for general $r \geq 1$ we establish that~:

\proclaim{Theorem 3.3} Suppose that $f$ satisfies {\rm ($H_1$)} with
$m \leq \frac{2r}{r+1} < q$, and $b$ satisfies {\rm ($H_2$)}-{\rm ($H_3$)}
for some $r\ge1$. 
If $\del = o(\eps^\frac{r+3}{r+1})$ then the sequence $u^{\eps,\del}$ 
converges in $L^s\left((0,T);L^p(\Rd) \right)$, for all $s<\infty$ and 
$p<q$, to a function
$$
u\in L^\infty\left((0,T);L^1(\Rd) \cap L^q(\Rd) \right),
$$  
which is the unique entropy solution to {\rm (1.1)}.
\endproclaim

Our method of proof can also be extended to a general diffusion 
$b=b(u,\grad u,D^2 u)$. 


\

\subheading{4. Convergence Proofs}

The superscripts $\eps$ and $\del$ are omitted in this section, 
except when necessary. 
In the proof, we make frequent use of the following computation. 
Multiply (1.4a) by $\eta'(u)$ where  $\eta : \RR \to \RR$ is a 
sufficiently smooth function and define
$q : \RR \to \Rd$ by $q_j' = \eta'\,f_j'$. We have
$$
\align
\pa_t \eta(u) = &  -\eta'(u) \, \div\!f(u)
				             + \eps \sum_j \pa_{x_j} \!\!
\left( \eta'(u)
\,\, b_j(\grad u) \right) -\pa_{x_j}\!\eta'(u) \,\, b_j(\grad u)
					           \\
                & + \del \sum_j \pa_{x_j} \!\! \left( \eta'(u)
\,\, \pa_{x_j}^2\!u \right) -\pa_{x_j}\!\eta'(u) \,\, \pa_{x_j}^2\!u
					           \\
				          = &  -\div q(u) + \eps \sum_j \pa_{x_j} \!\!
\left( \eta'(u) \,\,
b_j(\grad u) \right) -\eps \, \eta''(u) \sum_j \pa_{x_j}\!u \,\,
b_j(\grad u)
					           \\
					           & + \frac{\del}{2} \sum_j 2
\, \pa_{x_j}\!\!\left(
\eta'(u) \,\, \pa_{x_j}^2\!u \right) -\eta''(u) \,\, \pa_{x_j}\!\!\left(
\pa_{x_j}\!u \right)^2, 
\endalign
$$
thus
$$
\aligned 
\pa_t\eta(u)+\div q(u) =
\, & \eps \, \div\!\left( \eta'(u) \,\, b(\grad u) \right)
-\eps \, \eta''(u) \,\, \grad u \cdot\, b(\grad u) \\
    &  - \frac{\del}{2} \sum_j \eta''(u) \,\, \pa_{x_j}\!\!\left(
\pa_{x_j}\!u \right)^2 + \del \sum_j \pa_{x_j} \!\! \left( \eta'(u) \,\,
\pa_{x_j}^2\!u \right),
\endaligned 
\tag 4.1a
$$
The last two terms in the right hand side of (4.1a) 
take also the form 
$$
\frac{\del}{2} \sum_j \eta'''(u) \left( \pa_{x_j}\!u \right)^3
-3 \, \pa_{x_j}\!\!\left( \eta''(u) \left( \pa_{x_j}\!u \right)^2 \right)
+2 \, \pa_{x_j}^2\!\!\left( \eta'(u) \,\, \pa_{x_j}\!u \right).
\tag 4.1b
$$
When $\eta$ is convex, the term containing $\eta''(u)$ has a favorable sign~: 
the diffusion dissipates the entropy $\eta$. 

We begin by collecting fundamental energy estimates in several lemma. 

\proclaim{Lemma 4.1} Let $\alpha\geq 1$ be any real. 
Any solution of {\rm (1.4a)} satisfies, for $t\in [0,T]$, 
$$
\aligned 
& \int_{\Rd} \frac{|u(t)|^{\alpha+1}}{\alpha+1} \,dx
    + \alpha \, \eps \int_0^t\!\int_{\Rd} |u|^{\alpha-1} \,
    \grad u \cdot\, b(\grad u)\,dxds \\ 
& = \int_{\Rd} \frac{|u_0|^{\alpha+1}}{\alpha+1} \,dx     
-\frac{\alpha}{2} \, \del \int_0^t\!\int_{\Rd} |u|^{\alpha-1}
    \sum_j \pa_{x_j}\!\!\left( \pa_{x_j}\!u \right)^2 \,dxds. 
\endaligned 
\tag 4.2a
$$
For $\alpha \geq 2$, the last term in the above identity can be replaced by 
$$
\frac{\alpha\,(\alpha-1)}{2} \, \del
    \int_0^t\!\int_{\Rd} \sgn\!(u) \, |u|^{\alpha-2} \, \sum_j
    \left( \pa_{x_j}\!u \right)^3 \,dxds. 
\tag 4.2b
$$
\endproclaim

\proclaim{Proof} \rm Integrate (4.1a) over the whole of $\Rd$ with 
$\eta(u) = \frac{|u|^{\alpha+1}}{\alpha+1}$~:
$$
\aligned
& \frac{\t{d}}{\t{dt}} \int_{\Rd} \frac{|u|^{\alpha+1}}{\alpha+1} \,dx 
\\
& =
-\alpha \, \eps \int_{\Rd} |u|^{\alpha-1} \grad u \cdot b(\grad u)\,dx
-\frac{\alpha\, \del }{2} \int_{\Rd} \sum_j |u|^{\alpha-1}
\pa_{x_j}\!\!\left( \pa_{x_j}\!u \right)^2 \,dx,
\endaligned
$$
which yields (4.2a) after integration over $[0,t]$. 
One may use (4.1b), instead, to obtain (4.2b).
{$\square$}\medskip
\endproclaim

Choosing $\alpha = 1$ in Lemma 4.1, we deduce immediately 
a uniform bound for $u$ in
$L^\infty((0,T);L^2(\Rd))$ together with a control for both
$\grad u \cdot b(\grad u)$ in
$L^1(\Rd\times (0,T))$ and $\grad u$ in $L^{r+1}(\Rd\times  (0,T))$.

\proclaim{Proposition 4.2} For any solution of {\rm (1.4a)} and
$t \in [0,T]$, we have 
$$
\int_{\Rd} u(t)^2\,dx + 2\,\eps \int_0^t\!\int_{\Rd}
    \grad u \cdot b(\grad u)\,dxds
= \int_{\Rd} u_0^2\,dx 
\tag 4.3
$$
and, assuming {\rm ($H_2$)}, 
$$
\eps \int_0^t\!\int_{\Rd} |\grad u|^{r+1}\,dxds
\leq C \! \int_{\Rd} u_0^2\,dx.
\tag 4.4
$$
\endproclaim

To derive additional a~priori estimates, we use another value of $\alpha$, 
motivated by controling the dispersive term in (4.2b) with H\"older inequality, 
as follows~: 
$$
\aligned
& \left| \int_0^t\!\int_{\Rd} \sgn\!(u) \, |u|^{\alpha-2} \, \sum_j \left(
				\pa_{x_j}\!u \right)^3 \,dxds \right|
\\
&				\leq \int_0^t\!\int_{\Rd} |u|^{\alpha-2} \,
\left| \grad u \right|^3 \,dxds
                \\
&					 \leq \left[ \int_0^t\!\int_{\Rd}
|u|^{(\alpha-2)p} \,dxds
    \right]^\frac{1}{p} \left[ \int_0^t\!\int_{\Rd} \left|
    \grad u \right|^{3p'} dxds \right]^\frac{1}{p'}. 
\endaligned
\tag 4.5
$$
To take advantage of (4.4), we can choose $3 \, p'= r+1$ provided $r \geq 2$. 
Then $p= \frac{r+1}{r-2}$,
so $(\alpha-2)p= (r+1)\frac{\alpha-2}{r-2}$. 
Therefore it is rather natural to take the exponent $\alpha=r$ 
for the entropy, where $r$ is given by the diffusion term. 
Thus we deduce from Lemma 4.1 a natural estimate for $|u(t)|^{r+1}$,
involving the combination $\del \,\, \eps^{-\frac{3}{r+1}}$ of $\del$ and
$\eps$.

\proclaim{Proposition 4.3} Assume that {\rm ($H_2$)} holds with
$r\geq 2$ and $u_0\in L^{r+1}(\Rd)$. For $t \in [0,T]$, we have
$$
\aligned
& \int_{\Rd} |u(t)|^{r+1} \,dx + (r+1)\,r\,
    \eps \int_0^t\!\int_{\Rd} |u|^{r-1} \grad u \cdot b(\grad u)\,dxds \\
& \leq \,\, C_1(u_0) \left( 1 + \del \,\, \eps^{-\frac{3}{r+1}}\,
\max \left\{ 1, \left[ t \, C_1(u_0) 
\left( 1 + \del \,\, \eps^{-\frac{3}{r+1}} \right) \right]^\frac{r-2}{3}
\right\} \right) \\
& := \, H_{1} \! \left( \del \,\, \eps^{-\frac{3}{r+1}} \right)
\endaligned 
\tag 4.6 
$$
and 
$$
\eps \int_0^t\!\int_{\Rd}
    |u|^{r-1} \left| \grad u \right|^{r+1}\,dxds
    \leq  \frac {C}{(r+1) \, r}\, 
	 H_{1} \! \left( \del \,\eps^{-\frac{3}{r+1}}\right), 
\tag 4.7 
$$
where $C>0$ is some fixed constant and 
$$
C_1(u_0) := \max \left\{ \|u_0\|^{r+1}_{L^{r+1}(\Rd)},\,
         \frac{(r+1)r(r-1)}{2} \left( C \|u_0\|_{L^2(\Rd)}^2
			\right)^\frac{3}{r+1}\right\}.
$$ 
\endproclaim

In particular Proposition 4.3 shows that, if $u_0 \in L^2 \cap L^{r+1}$ 
and $\delta = O(\epsilon^{3 \over r+1})$, then $u(t) \in L^{r+1}$ 
uniformly for all $t \geq 0$. 

To motivate the forthcoming derivation, let us consider the 
special case $r=2$. Then (4.6) gives us an $L^3$ estimate. 
Returning to the original inequality (4.5), but now 
with the new value $\alpha = 3$, we now can estimate the dispersive 
term in (4.2b) directly in view of the estimate (4.7). 
In this fashion, we deduce an $L^4$ estimate from Lemma 4.1.
This argument can be continued inductively to reach any space $L^q$.

Actually Propositions 4.2 and 4.3 are the
first two cases of a general result derived now.
We define, for $n \geq 1$, 
$$ 
H_{0}\!\left( \del \,\, \eps^{-\frac{3}{r+1}} \right) =
   C_0(u_0) := \|u_0\|_{L^2(\Rd)}^2,
   $$
   $$
H_{n}\!\left( \del \,\, \eps^{-\frac{3}{r+1}} \right) :=
   C_n(u_0) \left( 1 + \del \,\, \eps^{-\frac{3}{r+1}} \,
   \max \left\{ 1,\left[ \, t \,\, C_n(u_0) 
   \left( 1 + \del \,\, \eps^{-\frac{3}{r+1}} \right)
	  \right]^\frac{r-2}{3} \right\} \right),
	  $$
	  and 
	  $$
\aligned
C_n(u_0):= 
&   \max \left\{ \|u_0\|^{n(r-1)+2}_{L^{n(r-1)+2}(\Rd)},\,\,
	  \frac{n(r-1)+2}{\left[ (n-1)(r-1)+2 \right]^\frac{3}{r+1}}
	  \frac{n(r-1)+1}{\left[ (n-1)(r-1)+1 \right]^\frac{3}{r+1}}
	  \right.
                 \\			
&			\left. n \, \frac{r-1}{2} \left(
			C \,H_{n-1} \! \left( \del\,\,\eps^{-\frac{3}{r+1}}
\right)
      	\right)^\frac{3}{r+1} \right\}. 
\endaligned
\tag 4.8
$$
Here $C>0$ is some fixed constant. Note that $H_n$ and $C_n$ are uniformly
bounded in $\epsilon, \delta$ provided $u_0 \in L^2\cap L^{n(r-1)+2}$ 
and $\delta= O(\epsilon^{3 \over r+1})$.

\proclaim{Proposition 4.4} Assume that {\rm ($H_2$)} holds with
$r\geq 2$ and $u_0\in L^q(\Rd)$. For $t \in [0,T]$ and $n \geq 0$
such that $n(r-1)+2 \leq q$, we have
$$
\aligned
& \int_{\Rd} |u(t)|^{n(r-1)+2} \,dx\\
& + \eps \, (n(r-1)+2)(n(r-1)+1) 
   \int_0^t\!\int_{\Rd} |u|^{n(r-1)} \,
    \grad u \cdot b(\grad u)\,dxds
\\
& \leq H_{n} \! \left( \del \,\, \eps^{-\frac{3}{r+1}} \right),
\endaligned
\tag 4.9
$$
and 
$$
\aligned
& \eps \int_0^t\!\int_{\Rd} |u|^{n(r-1)} \left| \grad u \right|^{r+1}\,dxds
\\
& \leq C \, (n(r-1)+2)^{-1}(n(r-1)+1)^{-1} H_{n} 
	 \! \left( \del \,\, \eps^{-\frac{3}{r+1}} \right). 
\endaligned
\tag 4.10
$$
\endproclaim

\proclaim{Proof of Propositions 4.3 and 4.4} 
\rm Note first that (4.10) is an immediate consequence
of (4.9) and the hypothesis ($H_2$). If  $n=0$, (4.9) coincides with 
(4.3) in Proposition 4.2. For $n=1$, the estimate is Proposition 4.3.  

To estimate the term in (4.2b), with $\alpha = r$, we use (4.5)~:
$$
\aligned
& \int_{\Rd} |u(t)|^{r+1}\,dx
+ (r+1)\,r\, \eps \int_0^t\!\int_{\Rd} |u|^{r-1}\,
    \grad u \cdot b(\grad u)\,dxds
    \\
& 	 \leq \int_{\Rd} |u_0|^{r+1}\,dx
               \\
& + \frac{(r+1)r(r-1)}{2} \, \del
   \left[ \int_0^t\!\int_{\Rd} |u|^{r+1}\,dxds \right]^\frac{r-2}{r+1}
   \left[ \int_0^t\!\int_{\Rd} \left| \grad u \right|^{r+1} \,dxds
   \right]^\frac{3}{r+1}.
\endaligned 
\tag 4.11
$$
By {\rm ($H_2$)} the second term in the left hand side of (4.11) 
is positive. Integrate (4.11) over $[0,t]$
and use (4.4)~:
$$
\align
& \|u\|^{r+1}_{L^{r+1}(\Rd\times (0,T))} 
\\
& \leq t \, \|u_0\|^{r+1}_{L^{r+1}(\Rd)}
               \\
& \quad + \, \frac{(r+1)r(r-1)}{2} \,\, t \, \left(
   C \, \|u_0\|^{2}_{L^2(\Rd)} \right)^{\frac{3}{r+1}} \,
	  \del \,\, \eps^{-\frac{3}{r+1}} \,\,
   \|u\|^{r-2}_{L^{r+1}(\Rd\times (0,T))}
                \\
					 & \leq  t \,\, C_1(u_0) \left( 1 +
\del \,\, \eps^{-\frac{3}{r+1}}
    \, \left( \|u\|^{r+1}_{L^{r+1}(\Rd\times (0,T))}
	   \right)^\frac{r-2}{r+1} \right). 
\endalign
$$

Observe that, for $X>0$, the inequality 
$$
X \leq \, K \left( 1 + \Delta \, X^\frac{\theta}{r+1}
    \right).
$$
where $0 \leq \theta < r+1$ and $K>0$, implies 
$$
X \leq \, \max \left\{ 1,\, \left[ \, K \left( 1 + \Delta
     \right) \right]^\frac{r+1}{r+1-\Theta} \right\}.
\tag 4.12
$$
Thus we deduce
$$
\|u\|^{r+1}_{L^{r+1}(\Rd\times (0,T))} \leq \, \max \left\{ 1,\,
    \left[ \, t \,\, C_1(u_0) \left( 1 + \del \, \eps^{-\frac{3}{r+1}}
    \right) \right]^\frac{r+1}{3} \right\} 
$$
and, returning to (4.11):
$$
\aligned
& \int_{\Rd} |u(t)|^{r+1}\,dx + (r+1)\,r\, \eps \int_0^t\!\int_{\Rd}
   |u|^{r-1} \grad u \cdot b(\grad u)\,dxds
                \\
&					 \leq C_1(u_0) \left( 1 + \del
\, \eps^{-\frac{3}{r+1}} \,
   \max \left\{ 1,\, \left[ \, t \,\, C_1(u_0) \left(
   1 + \del \, \eps^{-\frac{3}{r+1}} \right) \right]^\frac{r-2}{3} \right\}
   \right)
                \\
& := H_{1} \! \left( \del \,\, \eps^{-\frac{3}{r+1}} \right).
\endaligned
$$
This completes the proof of (4.6).

This argument can be iterated. 
We return to the dispersive term and make an estimate similar to (4.5),
but now having in view to apply (4.10), already established for $n=1$~:
$$
\aligned
& \left| \int_0^t\!\int_{\Rd}
\sgn\!(u) |u|^{\alpha-2}\, \sum_j \left(
					\pa_{x_j}\!u \right)^3 \,dxds \right|
\\
& \le
\int_0^t\!\int_{\Rd} |u|^{\alpha-2}\, \left| \grad u \right|^3 \,dxds
\\
& \leq \left[ \int_0^t\!\int_{\Rd} |u|^{(\alpha-2-\gamma)p} \,dxds
    \right]^\frac{1}{p} \left[ \int_0^t\!\int_{\Rd}
	 |u|^{\gamma p'} \left|
	 \grad u \right|^{3p'} \,dxds \right]^\frac{1}{p'},
\endaligned 
\tag 4.13
$$
where we choose $3p'=r+1$ and $\gamma p'=r-1$, so
$(\alpha-2-\gamma)p = \left( \alpha -2
        -3\,\frac{r-1}{r+1} \right) \frac{r+1}{r-2}\,$.
Then (4.2b) gives
$$
\aligned
& \int_{\Rd} |u(t) |^{\alpha+1}\,dx
+ (\alpha+1)\,\alpha\,\eps \int_0^t\!\int_{\Rd} |u|^{\alpha-1}\,
      \grad u \cdot b(\grad u)\,dxds
\\
&		\leq \int_{\Rd} |u_0|^{\alpha+1} dx
               \\
 & \quad  + \frac{(\alpha+1)\alpha(\alpha-1)}{2  \left[
  (r+1)r \right]^{\frac{3}{r+1}}} \left(
  CH_{1}\!\!\left( \del\,\,\eps^{-\frac{3}{r+1}} \right)
  \right)^\frac{3}{r+1} \del\,\eps^{-\frac{3}{r+1}}
  \left[ \int_0^t\!\int_{\Rd} |u|^{(\alpha-2-\gamma)p} dxdt
  \right]^\frac{r-2}{r+1}.
\endaligned
\tag 4.14
$$
We choose $\alpha$ so that $\alpha+1=(\alpha-2-\gamma)p\,$, i.e.,
$\alpha=2r-1\,$.

Integrating (4.14) over the interval $[0,t]$, we obtain
$$
\aligned 
& \|u\|^{2r}_{L^{2r}(\Rd\times (0,T))} 
\\
& \le \,\, t \,\, \|u_0\|^{2r}_{L^{2r}(\Rd)}
               \\
& \quad + \, \frac{r \, (2r-1)  \, (2r-2)}{\left[
   (r+1) \, r \right]^\frac{3}{r+1}} \,\, t \, \left(
   C \, H_{1} \!\! \left( \del \,\, \eps^{-\frac{3}{r+1}} \right)
   \right)^\frac{3}{r+1} \del \,\, \eps^{-\frac{3}{r+1}}
   \left( \|u\|^{2r}_{L^{2r}(\Rd\times (0,T))} \right)^\frac{r-2}{r+1} \\
&\leq 
 \,\, t \,\, C_2(u_0) \left( 1 + \del \,\, \eps^{-\frac{3}{r+1}}
    \, \left( \|u\|^{2r}_{L^{2r}(\Rd\times (0,T))}
	 \right)^\frac{r-2}{r+1} \right),
\endaligned 
$$
with 
$C_2(u_0) := \max \left\{ \|u_0\|^{2r}_{L^{2r}(\Rd)},\,
     \frac{r \, (2r-1) \, (2r-2)}{\left[
     (r+1) \, r \right]^\frac{3}{r+1}} \left(
     C \, H_{1} \!\! \left( \del \,\, \eps^{-\frac{3}{r+1}} \right)
     \right)^\frac{3}{r+1} \right\}$. 

By (4.12), we obtain again
$$
\|u\|^{2r}_{L^{2r}(\Rd\times (0,T))}
\leq \max \left\{ \, 1, \left[ \, t \,\, C_2(u_0)
    \left( 1 + \del \,\, \eps^{-\frac{3}{r+1}} \right)
	   \right]^\frac{r+1}{3} \right\}.
$$
Then (4.14) gives 
$$
\align
& \int_{\Rd} |u(t)|^{2r}\,dx
+ 2r \, (2r-1) \, \eps  \int_0^t\!\int_{\Rd} |u|^{2(r-1)} \,
    \grad u \cdot b(\grad u)\,dxds \\
&					  \leq  \, C_2(u_0) \left( 1 + \del
\,\, \eps^{-\frac{3}{r+1}}
    \max \left\{ 1, \left[ t \,\, C_2(u_0)
    \left( 1 + \del \,\, \eps^{-\frac{3}{r+1}} \right)
	   \right]^\frac{r-2}{3} \right\} \right)\\
& :=  \, H_{2} \! \left( \del \,\, \eps^{-\frac{3}{r+1}} \right).
\endalign
$$
This establishes (4.9) for $n=2$. The general case follows by induction on $n$. 
{$\square$}\medskip
\endproclaim

We are now concerned with the case where the diffusion
exponent in {\rm ($H_2$)} satisfies $r < 2$. In this situation, 
we require the assumption {\rm ($H_3$)}, which for instance is satisfied by 
$b_j(\grad u) = \partial_{x_j} u$.

\proclaim{Proposition 4.5} Suppose that 
{\rm ($H_1$)}--{\rm ($H_3$)} hold 
with $m$ and $r$ such that $m \leq \frac{2r}{r+1}$ and $r\ge 1$.
For $t\in [0,T]$, we have 
$$
\eps^\frac{r+3}{r+1} \int_{\Rd} \left| \grad u (t)\right|^2\,dx\,
+ \eps^{\frac{2(r+2)}{r+1}} \int_0^T\!\int_{\Rd}
\left|D^2 u\right|^2\,dxdt \,\le\, C, 
\tag 4.15
$$
$$
\int_{\Rd} \left| u (t) \right|^{2+\frac{r-1}{r}}\,dx\,
+ \eps \int_0^T\!\int_{\Rd} |u|^{\frac{r-1}{r}}\, |\grad u|^{r+1} \,dxdt 
  \leq \, C \left(\,1 + \, \del^\frac{r+1}{r} \, \eps^{-\frac{r+3}{r}}\right). 
\tag 4.16
$$
\endproclaim

\proclaim{Proof} \rm We differentiate (1.4a) with respect to the space variable $x$~: 
$$
\pa_{t} \grad u + \div \left( f'(u).\grad u \right) =
\eps \grad \sum_j \pa_{x_j} \left(b_j(\grad u) \right)
+ \del \sum_j \pa_{x_j}^3 \left( \grad u \right)
$$
and, then, we multiply it by $\grad u$ and integrate over $\Rd$. 
After further integration by parts, we obtain 
$$
\aligned 
& \frac{1}{2}\frac{\t{d}}{\t{dt}}\int_{\Rd} \left| \grad u (t)\right|^2\,dx
  -\int_{\Rd} \lap u \,f'(u) \cdot \grad u \,dx \\
& = 
  -\eps \int_{\Rd} \sum_k \grad \partial_{x_k} u \cdot Db(\grad u) \cdot 
  \grad \partial_{x_k} u \,dx
-\frac{\del}{2} \sum_j \int_{\Rd} \pa_{x_j}\!\left(
\sum_k \left( \pa_{x_k x_j}^2 u\right)^2 \right).
\endaligned 
$$
Thus, integrating on $[0,t]$ using ($H_1$) yields 
$$
\aligned 
& \int_{\Rd} \left| \grad u (t) \right|^2\,dx \, 
  + 2 \, \eps \int_{\Rd} \sum_k \grad \partial_{x_k} u \cdot Db(\grad u) \cdot 
  \grad \partial_{x_k} u \,dx\\ 
& \leq 
   \int_{\Rd} \left| \grad u_0 \right|^2\,dx \, 
   + \, 2 \, C_1\int_0^t\!\int_{\Rd}
    \left| D^2 u \right| \left| u \right|^{m-1} 
    \left|\grad u\right| \,dxdt \\ 
&   \leq 
    \int_{\Rd} \left| \grad u_0 \right|^2\,dx \,
   + \, \frac{C}{\eps}\int_0^t\!\int_{\Rd}
    \left| u \right|^{2m-2} \left|\grad u\right|^2\,dxdt \,
   + \, C_4 \, \eps\int_0^t\!\int_{\Rd}
\left|D^2 u\right|^2\,dxdt,
\endaligned 
$$
and so, using ($H_3$),
$$
\aligned 
\int_{\Rd} \left| \grad u (t) \right|^2\,dx \,
& + C_5 \, \eps \int_0^t\!\int_{\Rd}
\left|D^2 u\right|^2\,dxdt \\
 & \le \int_{\Rd} \left| \grad u_0 \right|^2\,dx \,
+ \, \frac{C}{\eps}\int_0^t\!\int_{\Rd}
      \left| u \right|^{2m-2} \left|\grad u\right|^2\,dxdt. 
\endaligned 
$$
By H\"older inequality and for $m \leq \frac{r-1}{r+1}$~:
$$
\multline
\int_{\Rd} \left| \grad u (t)\right|^2\,dx \,
+  C_5 \, \eps\int_0^t\!\int_{\Rd}
     \left| D^2 u\right|^2\,dxdt
	  \\
	  \leq \int_{\Rd} \left| \grad u_0 \right|^2\,dx \,
+ \, C \, \eps^{-1} \left[ \int_0^t\!\int_{\Rd}
      \left|\grad u\right|^{r+1}\,dxdt \right]^\frac{2}{r+1}
\left[ \int_0^t\!\int_{\Rd}
      \left| u \right|^2 \,dxdt \right]^\frac{r-1}{r+1},
\endmultline
$$
and now (4.15) follows from (4.3)-(4.4).

To establish (4.16)-(4.17) we use (4.2a) for $\alpha \ge 1$~:
$$
\aligned 
& \int_{\Rd} |u(t)|^{\alpha+1}\,dx \,
    +  C \,\eps \int_0^t\!\int_{\Rd} |u|^{\alpha-1}\, |\grad u)|^{r+1} \,dxdt
   \\
& \leq \int_{\Rd} |u_0|^{\alpha+1}\,dx \,
    +C' \, \del \int_0^t\!\int_{\Rd} |u|^{\alpha-1} \, \grad u| \, |D^2 u| \,dxdt. 
\endaligned 
$$
We evaluate the last term using ($H_2$):
$$
\aligned 
&  \del\int_0^t\!\int_{\Rd} |u|^{\alpha-1} \, \grad u| \, |D^2 u| \,dxdt
   \\
&  \leq \del \int_0^t\!\int_{\Rd} |u|^{\alpha-1} \left(
      \frac{C_2 \, \eps}{(r+1) \, \del} \left| \grad u \right|^{r+1}\,
		+ \, \frac{r}{r+1} \left( \frac{\del}{C_2 \, \eps}
\right)^\frac{1}{r}
		\left| D^2 u \right|^\frac{r+1}{r} \right) \,dxdt \\
&	\leq \frac{\eps}{2} \int_0^t\!\int_{\Rd}
    |u|^{\alpha-1} \, \left| \grad u \right|^{r+1}\,dxdt \,
	 + C'' \, \del^\frac{r+1}{r}\eps^{-\frac{1}{r}}
	 \int_0^t\!\int_{\Rd} |u|^{\alpha-1}
	 \left|D^2 u \right|^\frac{r+1}{r}\,dxdt. 
\endaligned 
$$
So, we have 
$$
\aligned 
& \int_{\Rd} |u(t)|^{\alpha+1}\,dx +  C \, \eps
    \int_0^t\!\int_{\Rd} |u|^{\alpha-1}\, |\grad u|^{r+1} \,dxdt
   \\
	& \le \int_{\Rd} |u_0|^{\alpha+1}\,dx + C \,
    \del^\frac{r+1}{r}\eps^{-\frac{1}{r}}
	   \int_0^t\!\int_{\Rd} |u|^{\alpha-1}
	   \left|D^2 u \right|^\frac{r+1}{r}\,dxdt. 
\endaligned 
$$
Taking $\alpha = 1+\frac{r-1}{r}$, we deduce
$$
\aligned 
& \int_{\Rd} |u(t)|^{2+\frac{r-1}{r}}\,dx + C \, \eps
    \int_0^t\!\int_{\Rd} |u|^{\frac{r-1}{r}}\, |\grad u|^{r+1} \,dxdt \\ 
& \le \int_{\Rd} |u_0|^{2+\frac{r-1}{r}}\,dx + C \, 
    \del^\frac{r+1}{r}\eps^{-\frac{1}{r}}
	   \left[ \int_0^t\!\int_{\Rd} |u|^{2}\,dxdt \right]^\frac{r-1}{2r}
	   \left[ \int_0^t\!\int_{\Rd} \left|D^2 u
    \right|^2\,dxdt \right]^\frac{r+1}{2r}.
\endaligned 
$$
The conclusion follows now easily. 
{$\square$}\medskip
\endproclaim


\proclaim{Proof of Theorems 3.1} \rm We first derive (2.5a), based on 
the conservation law (4.1b) with
a arbitrary convex function, $\eta$, where we assume
$\eta ',\eta '',\eta '''$ bounded functions on $\RR$.
We claim that there exists a bounded 
measure $\mu \leq 0$ such that
$$
\pa_t \eta(u) + \div q(u) \longrightarrow
       \mu \quad \t{ in } \Cal D'(\Rd\times (0,T)).
$$
{}From (4.1b), we obtain  
$$
\split
\pa_t\eta(u)+\div q(u) = 
   &  \,\, \eps \, \div\!\left( \eta'(u) \,\,
      b(\grad u) \right) -\eps \, \eta''(u) \,\,
	     \grad u \cdot b(\grad u)
								        \\
   &  +\frac{\del}{2} \sum_j \eta'''(u) \, \left( \pa_{x_j}\!u \right)^3 
      -3 \, \pa_{x_j} \!\! \left( \eta''(u) \,
      \left( \pa_{x_j}\!u \right)^2 \right)
      +2\,\pa_{x_j}^2 \!\! \left( \eta'(u) \,\,
      \pa_{x_j} \! u \right) \\
:= &  \, \mu_1 + \mu_2 + \mu_3,
\endsplit
$$
with obvious notation. For each positive
$\theta \in \CC^\infty_0(\Rd\times (0,T))$
we evaluate $\la \mu_i,\theta \ra$ for $i = 1,2,3$.
To treat $\mu_1$, we use H\"older inequality with the exponent 
$\frac{r+1}{r}$. In view of ($H_2$) and 
(4.4) of Proposition 4.2 and assumption (3.1), we get 
$$
\align
\left| \la \mu_1,\theta \ra \right|
   &  \leq \eps \int^T_0\!\!\int_{\Rd} \sum_j \left|
	     \eta'(u)\, b_j(\grad u)
		  \,\,\pa_{x_j}\!\theta \right|\,dxdt
\\
   &  \leq C\,\eps \int^T_0\!\!\int_{\Rd}
		  \left| \grad\theta \right|\,\left|
		  b(\grad u) \right|\,\,dxdt
\endalign
$$
so 
$$
\align 
\left| \la \mu_1,\theta \ra \right| 
   &  \leq C\,\eps\, \|\grad\theta\|_{L^{r+1}(\Rd\times (0,T))}
		  \left[ \iint_{\supp\theta} \left| \grad u \right|^{r+1}
		  \,dxdt \right]^\frac{r}{r+1} \\
	&  \leq C\,\eps^\frac{1}{r+1}\,
	  \|\grad\theta\|_{L^{r+1}(\Rd\times (0,T))}.
\endalign
$$
For $\mu_2$, we use ($H_2$) and the convexity of $\eta$~: 
$$
\la \mu_2,\theta \ra = - \eps \int^T_0\!\!\int_{\Rd} \sum_j
                         \theta \,\, \eta''(u) \,\,
                         \grad u \cdot b(\grad u)  \,\,dxdt 
								 \leq 0.
$$
For $\mu_3$, we use again H\"older inequality, as follows  
$$
\align
& \left| \la \mu_3,\theta \ra \right|
\\  & \leq \,\frac{\del}{2} \int^T_0\!\!\int_{\Rd} \sum_j 
		  \left| \theta\,\,\eta'''(u)\,\left( \pa_{x_j}\!u \right)^3
		  +3\,\eta''(u) \left( \pa_{x_j}\!u \right)^2
		  \pa_{x_j}\!\theta
		  + \,2\,\eta'(u)\,\,\pa_{x_j}\!u\,
		  \,\pa_{x_j}^2\!\theta \right| \,dxdt
\\
  & \leq \,C\,\del \int^T_0\!\!\int_{\Rd} \theta\,\left|
	     \grad u \right|^3 \,dxdt
	     +\,C\,\del \int^T_0\!\!\int_{\Rd} \sum_j \left|
		  \pa_{x_j}\!u \right|^2 \left| \pa_{x_j}\!\theta \right|
		  \,dxdt
\\
  & \quad +\,C\,\del \int^T_0\!\!\int_{\Rd} \sum_j \left|
		  \pa_{x_j}\!u \right| \,\left| \pa_{x_j}^2 \!\theta \right|
		  \,dxdt, 
\endalign
$$
so 
$$
\align
\left| \la \mu_3,\theta \ra \right|
  & \leq \,\,C\,\del \,
        \| \theta\|_{L^\frac{r+1}{r-2}(\Rd\times (0,T))}\,\,\left[
		  \iint_{\supp\theta}\left| \grad u \right|^{r+1} \,dxdt
		  \right]^\frac{3}{r+1}
\\
  & \quad +\,C\,\del \,
        \| \grad\theta\|_{L^\frac{r+1}{r-1}(\Rd\times (0,T))}
		  \,\,\left[
		  \iint_{\supp\theta}\left| \grad u \right|^{r+1} \,dxdt
		  \right]^\frac{2}{r+1}
\\
  & \quad +\,C\,\del \,\left[ \int^T_0\!\!\int_{\Rd} \left( \sum_j \left|
		  \pa_{x_j}^2\!\theta \right| \right)^\frac{r+1}{r} dxdt
		  \right]^\frac{r}{r+1}
		  \,\,\left[
		  \iint_{\supp\theta}\left| \grad u \right|^{r+1} \,dxdt
		  \right]^\frac{1}{r+1}. 
\endalign
$$
Therefore, we conclude that  
$$
\left| \la \mu_3,\theta \ra \right|
 \leq \,\,C\,\del \left( \eps^{-\frac{3}{r+1}}
		  +\,\eps^{-\frac{2}{r+1}}
		  +\,\eps^{-\frac{1}{r+1}} \right)
		  \leq \,C\,\del\,\eps^{-\frac{3}{r+1}}.
$$
Finally the condition $\del = o(\eps^\frac{3}{r+1})$ is sufficient to 
imply the desired conclusion.

Using a standard regularization of $\sgn\!(u)$ and $|u-k|$
(for $k\in\RR$), which fullfil the growth condition (2.1), we apply
the limit representation (2.2) and conclude that $\nu$ satisfies (2.5a).

To show (2.5b) we follow DiPerna \cite{\refDiPerna}
and Szepessy \cite{\refSzepessy}'s arguments. 
We have to check that, for each compact $K$ of $\Rd$,
$$
\aligned 
& \lim_{t\to 0+} \frac{1}{t} \int_0^t\!\int_K
               \la \nu_{(x,s)}, \left| u-u_0(x) \right|\,\ra\,\,dxds \\ 
& = \lim_{t\to 0+} \lim_{\eps\to 0+}\frac{1}{t}
	\int_0^t\!\int_K
   \left| u^{\eps,\del}(x,s)-u_0(x) \right| \,dxds
	= 0.
\endaligned 
$$
By Jensen's inequality, where $m(K)$ stands for Lebesgue
measure of $K$, we have 
$$
\aligned
& \frac{1}{t} \int_0^t\!\int_K
     \left| u^{\eps,\del}(x,s)-u_0(x) \right| \,dxds
\\
&	  \leq m(K)^{1/2} \left( \frac{1}{t} \int_0^t\!\int_K
     \left( u^{\eps,\del}(x,s)-u_0(x) \right)^2 \,dxds
	  \right)^{1/2}. 
\endaligned
$$
We will establish that
$$
\lim_{t\to 0+} \lim_{\eps\to 0+}
\frac{1}{t} \int_0^t\!\int_K
     \left( u^{\eps,\del}(x,s)-u_0(x) \right)^2 \,dxds = 0.
$$
Let $K_{i} \subset K_{i+1}$ (\,$i = 0,1,...$\,) be an increasing
sequence of compact sets such that $K_0 = K$ and
$\cup_{i \ge 0} K_i = \Rd$. We use the identity
$u^2 - u_0^2 - 2u_0(u - u_0) = (u - u_0)^2\,$:
$$
\aligned 
& \frac{1}{t} \int_0^t \!\int_K
     \left( u^{\eps,\del}(\cdot,s)-u_0 \right)^2 \,\,dxds \\
& \leq \frac{1}{t} \int_0^t \left( \int_{K_i}
     |u^{\eps,\del}(\cdot,s)|^2 \,dx - \int_{K_i} u_0^2\,\,dx
	  - 2\int_{K_i} u_0 \left(u^{\eps,\del}(\cdot,s) - u_0 \right)
	  \,dx \right)\,ds \\
& \leq \int_{\Rd\setminus K_i} u_0^2\,dx
     + \frac{2}{t} \int_0^t \left| \int_{K_i}
	   u_0 \left(u^{\eps,\del}(\cdot,s) - u_0 \right)
	  \,\,dx \right|\,ds
\endaligned 
$$
for all $i = 0,1,...$, where we used (4.3)-(3.1).

Since 
$$
\lim_{i\to\infty} \int_{\Rd\setminus K_i} u_0^2\,\,dx = 0,
$$
we only consider the last term above. Take
$\{ \theta_n \}_{n\in\NN}\subset \CC^\infty_0(\Rd)$ such that
$$
\lim_{n\to\infty} \theta_n = u_0 \quad \t{ in } L^2(\Rd),
$$
Cauchy-Schwarz inequality gives
$$
\align
& \left| \int_{K_i} u_0 \left(u^{\eps,\del}(\cdot,s) - u_0 \right)
	  \,\,dx \right|
\\
&
\leq \int_{K_i} \left| u_0 - \theta_n \right|
     \left| u^{\eps,\del}(\cdot,s) - u_0 \right| \,\,dx
\\
&
\quad + \left| \int_{K_i} \theta_n \left( u^{\eps,\del}_0 - u_0 \right)
	  + \int_{K_i} \theta_n \left( u^{\eps,\del}(\cdot,s)
	  - u^{\eps,\del}_0 \right) \,dx \right|
\\
&
\leq \| u_0 - \theta_n \|_{L^2(\Rd)} \left(
	  \| u^{\eps,\del}(\cdot,s) \|_{L^2(\Rd)}
	  + \| u_0 \|_{L^2(\Rd)} \right)
\\
&
\quad + \| \theta_n\|_{L^2(\Rd)} \|u^{\eps,\del}_0 - u_0 \|_{L^2(\Rd)}
     + \left| \int_0^s\!\int_{K_i} \theta_n\,\pa_s u^{\eps,\del}
	  \,dxd\tau \right|.
\endalign
$$
In view of (4.3) and (3.1)
$$
\aligned 
& \| u_0 - \theta_n \|_{L^2(\Rd)} \left(
	  \| u^{\eps,\del}(\cdot,s) \|_{L^2(\Rd)}
	  + \| u_0 \|_{L^2(\Rd)} \right)
\\	  & \leq 2 \| u_0 \|_{L^2(\Rd)} \| u_0 - \theta_n \|_{L^2(\Rd)}, 
\endaligned 
$$
which tends to zero when $n \to \infty$, 
and since 
$\lim_{\eps\to 0+} \| u_0^{\eps,\del} - u_0 \|_{L^2(\Rd)} = 0$
by (3.1), it remains only to see that
$$
\lim_{t\to 0+}\,\lim_{\eps\to 0+}\frac{1}{t} \int_0^t \left|
    \int_0^s\!\int_{K_i} \theta_n\,\pa_s u^{\eps,\del}
	  \,dxd\tau \right| = 0.
$$
We have, by (1.4a), 
$$
\aligned 
\left| \int_0^s\!\int_{K_i} \theta_n\,\pa_s u\,dxd\tau \right|
& = \left| \int_0^s\!\int_{K_i} \theta_n \big( -\div\!f(u)
   + \eps\,\div b
	- \del\,\sum_j \pa_{x_j}^3 u \big)\,dxd\tau \right|
\\
& = \left| \int_0^s\!\int_{K_i} \big( \grad\theta_n \cdot f(u)
    - \eps\,\grad\theta_n \cdot b 
	 + \del \sum_j \pa_{x_j}^3 \theta_n \big) u 
	 \,dxd\tau \right|
\\
& := \mu_1 + \mu_2 + \mu_3. 
\endaligned 
$$
To deal with $\mu_1$, we use H\"older inequality and ($H_1$)
$$
\aligned 
& \int_0^s\!\int_{K_i} \left| \grad\theta_n \right|\,\left| f(u) \right|
\,dx\,d\tau\\
& \leq C \int_0^s\!\int_{K_i} \left| \grad\theta_n \right| \,dx\,d\tau
+ C \left[ \int_0^s\!\int_{K_i}
      \left| \grad\theta_n \right|^\frac{q}{q-m} \,dx\,d\tau
		\right]^\frac{q-m}{q}
		\left[ \int_0^s\!\int_{K_i} |u|^q\,dx\,d\tau
\right]^\frac{m}{q}\\
& \leq C\,s\, \| \grad\theta_n\|_{L^1(\Rd)} + C\,s^\frac{q}{q-m}\,
      \| \grad\theta_n\|_{L^\frac{q}{q-m}(\Rd)}\,
		\| u\|^{m}_{L^q(\Rd\times (0,T))}. 
\endaligned 
$$
For $\mu_2$, using ($H_2$) and once more H\"older inequality
with(4.4)-(3.1), we get 
$$
\align
&\eps\,\int_0^s\!\int_{K_i} \left| \grad\theta_n \right|\,|b|
    \,\,dx\,d\tau
\\
& \leq C_3\,\eps\,\int_0^s\!\int_{K_i} \left| \grad\theta_n \right|
    \left| \grad u\right|^r\,dx\,d\tau
\\
& \leq C_3\,\eps\,\left[ \int_0^s\!\int_{K_i}
     \left| \grad\theta_n \right|^{r+1}\,dxd\tau \right]^\frac{1}{r+1}
     \left[ \int_0^s\!\int_{K_i}
     \left| \grad u \right|^{r+1}\,dxd\tau \right]^\frac{r}{r+1}
\\
& \leq C\,\eps^{1-\frac{r}{r+1}}\,s^\frac{1}{r+1}\,
     \| \grad\theta_n \|_{L^{r+1}(\Rd)}. 
\endalign
$$
Finally, for $\mu_3$, we use Cauchy-Schwarz inequality with (4.3)-(3.1)~:
$$
\align
& \del\int_0^s\!\int_{K_i} \left| u \sum_j \pa_{x_j}^3 \theta_n \right|
    \,dx\,d\tau
\\
& \leq \del \left[ \int_0^s\!\int_{K_i}
    \left| u \right|^2 \,dxd\tau \right]^\frac{1}{2}
    \left[ \int_0^s\!\int_{K_i}
    \left| \sum_j \pa_{x_j}^3 \theta_n \right|^2 dx\,d\tau
	 \right]^\frac{1}{2}
\\
& \leq \del\,s\, \| \grad^3\theta_n \|_{L^2(\Rd)}\,\| u_0\|_{L^2(\Rd)},
\endalign
$$
thus 
$$
\aligned
& \lim_{\eps\to 0+} \frac{1}{t} \int_0^t \left| \int_0^s\!\int_{K_i}
     \theta_n\,\pa_s u^{\eps,\del}\,dxd\tau \right| ds  =: X 
	\\
	&  \leq \lim_{\eps\to 0+}\frac{1}{t} \left( \frac{C}{2}\,t^2\,
    \| \grad\theta_n \|_{L^1(\Rd)} \right.
\\
&  \quad + C \left( \frac{q}{q-m}+1 \right)^{-1}t^{\frac{q}{q-m}+1}\,
    \| \grad\theta_n \|_{L^\frac{q}{q-m}(\Rd)}\,
	 \| u^{\eps,\del}\|^{m}_{L^q(\Rd\times (0,T))}
\\
&  \quad + \left. C \left( \frac{1}{r+1}+1 \right)^{-1}\!
    t^{\frac{1}{r+1}+1} \eps^\frac{1}{r+1}
    \| \grad\theta_n \|_{L^{r+1}(\Rd)}
	 + \frac{\del t^2}{2} \| \grad^3\theta_n \|_{L^2(\Rd)}\,
	 \| u_0 \|_{L^2(\Rd)} \right)
\endaligned 
$$
and so 
$$
\aligned
X 
\leq 
& C_n\left( t + t^\frac{q}{q-m} \lim_{\eps\to 0+}
    \| u^{\eps,\del}\|^{m}_{L^q(\Rd\times (0,T))} \right)
\\
\leq & C_n\,t + C_{\eps,\del} \,t^\frac{q}{q-m},
\endaligned
$$
where we have used (4.9) in Proposition 4.4. 
The desired conclusion when $t\to 0+$ follows.
{$\square$}\medskip
\endproclaim


\proclaim{Proof of Theorems 3.2 and 3.3} \rm 
In the previous proof, to establish (2.5a) we started with the identity 
(4.1) and the condition $\del = o(\eps^\frac{3}{r+1})$
as required, inmparticular to control the term in (4.1b). 
We now keep the form (4.1a) instead (4.1b)~: the terms
t$\mu_1$ and $\mu_2$ introduced in the previous proof do not change. 
We only need discuss $\mu_3$. It has now the form:
$$
- \del\sum_j \frac{\eta''(u)}{2} \pa_{x_j}(\pa_{x_j}\!u)^2
+ \del \sum_j \pa_{x_j}(\eta'(u) \pa_{x_j}^2\!u).
$$
The first term is bounded as follows 
$$
\align
\left| \del \int_0^T\!\int_{\Rd} \sum_j
     \theta \, \eta''(u) \pa_{x_j}\!u \, \pa_{x_j}^2\!u \,dxdt \right| \,
	  \le \, & \,C \,\del \int_0^T\!\int_{\Rd}
     \theta \left| \grad u \right|
	  \left| D^2 u\right| \,dxdt
	 \\
	 \le \, & \, C \,\del \int_0^T\!\int_{\Rd}
     \mu \left| D^2 u\right|^2
	  + \frac{1}{\mu} \left( \theta \left| \grad u \right| \right)^2
	   \,dxdt
	 \\
	 \le \, & \, C \,\del \left( \mu \, \eps^{-2\frac{r+2}{r+1}} \,
+ \, \frac{1}{\mu} \, \eps^{-\frac{2}{r+1}} \right)
\endalign
$$
using (4.15) and (4.4), and we take $\mu = \eps$ and 
$\del = o(\eps^\frac{r+3}{r+1})$. 

The second term in $\mu_3$ behaves better:
$$
\align
& \del \left| \int_0^T\!\int_{\Rd} \theta \sum_j
     \pa_{x_j}(\eta'(u) \pa_{x_j}^2\!u) \right.\left.\,dxdt \right| =: Y
     \\ 
& \leq  \,\del \left| \int_0^T\!\int_{\Rd} \sum_j
     \pa_{x_j}^2\!\theta\,\eta'(u) \pa_{x_j}\!u \,dxdt \right|
	 + \del \left| \int_0^T\!\int_{\Rd} \sum_j
     \pa_{x_j}\!\theta\,\eta''(u) \left( \pa_{x_j}\!u \right)^2
	  \,dxdt \right|
	 \\
	& \le  \,C \, \del\int_0^T\!\int_{\Rd} \left| \grad u \right|
     \left| D^2\theta \right| \,dxdt
	  + \, C \, \del \int_0^T\!\int_{\Rd} \left| \grad u \right|^2
     \left| \grad \theta \right| \,dxdt
	 \endalign
	 $$
	 thus 
	 $$
	 \aligned
Y 	 \le & \,C \, \del \left[ \int_0^T\!\int_{\Rd} \left| \grad u
\right|^{r + 1}
    \,dxdt \right]^\frac{1}{r + 1}
	  + \, C \, \del \left[ \int_0^T\!\int_{\Rd} \left| \grad u
\right|^{r+1}
     \,dxdt \right]^\frac{2}{r+1}
	 \\
	 \le & \,C \, \del \,\eps^{-\frac{2}{r+1}}.
\endaligned
$$
This completes the proof of Theorems 3.2 and 3.3. 
{$\square$}\medskip
\endproclaim 

\

\subheading{Acknowledgments} 
The authors were partially supported by the
French Ambassy in Portugal and by Junta Nacional de 
Investiga\c c\~ ao Cient\'\i fica e Tecnologica, Portugal. 
J.M.C. was also partially supported by the
Funda\c c\~ ao Calouste Gulbenkian. 
P.G.L. was partially supported by the Centre National de la Recherche Scientifique, 
and by a Faculty Early Career Development award (CAREER) 
from the National Science Foundation, under grants DMS 94-01003 and DMS 95-02766. 


\

\subheading{References}

\

\ref\key{\refCLF} 
\by \, \,  F. Coquel and P.G. LeFloch
\paper Convergence of finite difference schemes
for scalar conservation laws in several space dimensions: the
corrected antidiffusive flux approach
\jour Math. of Comp. \vol \rm  57 \yr 1991 \pages 169--210
\endref

\ref\key{\refDiPerna}
\by \, \,   R.J. DiPerna
\paper Measure-valued solutions to conservation laws
\jour Arch. Rational Mech. Anal.
\vol \rm  88
\yr 1985
\pages 223--270
\endref

\ref\key{\refHLone} 
\by \, \,  B.T. Hayes and P.G. LeFloch
\paper Nonclassical shock waves and kinetic relations~:
scalar conservation laws 
\jour Arch. Rational Mech. Anal. \vol\rm 139  
\yr 1997 \pages 1--56 
\endref

\ref\key{\refHLtwo} 
\by \, \,  B.T. Hayes and P.G. LeFloch
\paper Nonclassical shock waves and kinetic relations~:
finite difference schemes 
\jour SIAM Numer. Anal. \vol\rm 35
 \yr 1998
\pages 2169--2194  
\endref

\ref\key{\refKruzkov}
\by \, \,  S.N. Kru\v zkov
\paper First order quasilinear equations in several independent
variables
\jour Mat. Sb.
\vol \rm  81
\yr 1970
\pages 285--355; {\it Math. USSR Sb.} {10} (1970), 217--243
\endref

\ref\key{\refLL} 
\by \, \,   P.D. Lax and C.D. Levermore
\paper The small dispersion limit of the Korteweg-de Vries equation
\jour  Comm. Pure Appl. Math. \vol \rm  36\yr 1983
\pages 253--290 
\endref

\ref\key{\refLN}
\by \, \,   P.G. LeFloch and R. Natalini
\paper Conservation laws with vanishing nonlinear
diffusion and dispersion
\jour J. of Nonlinear Analysis \yr 1997 
\endref

\ref\key{\refVonNeumann}
\by \, \,   J. Von Neumann and R.D. Richtmyer 
\paper A method for the
numerical calculation of hydrodynamical shocks
\jour J. Appl. Phys. 
\vol \rm  21
\yr 1950
\pages 380--385
\endref

\ref\key{\refSchonbek}
\by \, \,   M.E. Schonbek
\paper Convergence of solutions to nonlinear dispersive equations
\jour Comm. Part. Diff. Equa.
\vol \rm  7 
\yr 1982
\pages 959--1000
\endref

\ref\key{\refSzepessy} 
\by \, \, A. Szepessy
\paper An existence result for scalar conservation laws using
measure-valued
solutions
\jour Comm. Part. Diff. Equa.
\vol \rm  14
\yr 1989
\pages 1329--1350
\endref

\ref\key{\refTartar} 
\by \, \,  L. Tartar
\paper The compensated compactness method applied to systems of
conservation laws
\inbook Systems of Nonlinear Partial
Differential Equations, J. M. Ball ed.,
NATO ASI Series, C. Reidel publishing Col., \yr 1983
\pages 263--285
\endref

\ref\key{\refVolpert} 
\by \, \,   A.I. Volpert 
\paper The space BV and quasilinear equations 
\jour Math. USSR Sb.
\vol \rm  2
\yr 1967
\pages 257--267
\endref 

\enddocument